\documentclass[12pt]{amsart}
\usepackage[english]{babel}
\usepackage[utf8]{inputenc}
\usepackage{amssymb,amsmath,amsfonts,amsthm,amscd,latexsym,graphicx,epsfig,colordvi}

\title{Semiprime Novikov algebras}
\author{A.S.Panasenko}
\address{Alexander Sergeevich Panasenko
\newline\hphantom{iii} Sobolev Institute of Mathematics,
\newline\hphantom{iii} pr. Koptyuga, 4,
\newline\hphantom{iii} 630090, Novosibirsk, Russia.
\newline\hphantom{iii} Novosibirsk State University,
\newline\hphantom{iii} Universitetskiy pr., 1,
\newline\hphantom{iii} 630090, Novosibirsk, Russia}
\email{a.panasenko@g.nsu.ru}%
\thanks{\sc Panasenko, A.S.,
Semiprime Novikov algebras}
\thanks{\copyright \ 2022 Panasenko A.S}
\thanks{\rm Acknowledgments: The reported study was funded by RFBR, project number 21-11-00286.}
\thanks{\it Received  June, 1, 2022, published  December, 1,  2022.}%

\begin{document}

\maketitle

\vspace{1cm}
\maketitle {\small
\begin{quote}
\noindent{\sc Abstract. } We study prime and semiprime Novikov algebras. We prove that prime nonassociative Novikov algebra has zero nucleus and center. It is well known that an ideal of an alternative (semi)prime algebra is (semi)prime algebra. We proved this statement for Novikov algebras. It implies that a Baer radical exists in a class of Novikov algebras. Also, we proved that a minimal ideal of Novikov algebra is either trivial, or a simple algebra. 
\medskip

 \end{quote}
}

\section{Introduction.}

Simple and prime algebras are the main part of the structure theory of any variety. Novikov algebras first appeared in work \cite{GD} by I.~M.~Gelfand and I.~Ya.~Dorfman. Simple finite dimensional Novikov algebras over a field of characteristic zero were described by E.~I.~Zelmanov in \cite{Zelm}: they all are fields. Simple infinite dimensional Novikov algebras were studied in \cite{Osborn2},\cite{Xu1}. Simple finite dimensional Novikov algebras with an idempotent over a perfect field of characteristic $p>2$ were described in \cite{Osborn1}. The final description of simple finite dimensional Novikov algebras over an algebraically closed field of characteristic $p>2$ is given in \cite{Xu2}. 

As we know, prime and semiprime Novikov algebras has never been studied. Prime alternative algebras over a field of characteristic $\neq 3$ were described by Slater in \cite{Slater}. Prime nondegenerate Jordan algebras were described by E.I.Zelmanov in \cite{Zelm2}.

Novikov algebras were actively studied recently. I.~P.~Shestakov and Z.~Zhang proved in \cite{ShestZhang} the solvability and right-nilpotency are equivalent in the variety of Novikov algebras. Moreover, it is equivalent to nilpotency of the square of an algebra. U.~U.~Umirbaev and V.~N.~Zhelyabin \cite{UZ} proved an analogue of Bergman-Isaacs theorem: $Z_n$-graded Novikov algebra with a solvable $0$-component is solvable. K.M.Tulenbaev, U.U.Umirbaev and V.N.Zhelyabin \cite{TUZ} proved that a Novikov algebra over a field of characteristic $\neq 2$ is solvable iff its commutator ideal is right-nilpotent. 

Classical examples of prime algebras in varieties of associative, alternative and Jordan algebras are central orders in a simple unital algebra. But unital Novikov algebra is associative and commutative. So, we need to study some properties of prime Novikov algebras to understand, can we create a method to construct prime Novikov algebras.

In this paper, we proved that prime nonassociative Novikov algebra has zero nucleus and center. We proved that an ideal of (semi)prime Novikov algebra is (semi)prime Novikov algebra. It implies that a Baer radical exists in a class of Novikov algebras. Besides, it has been proven that a minimal ideal of Novikov algebra is either trivial, or a simple algebra. 

\section{Nucleus of prime Novikov algebras.}

We will use usual notations for associator and commutator:
\[(a,b,c)=(ab)c-a(bc), \qquad [a,b]=ab-ba.\]

Notation $I\unlhd A$ means that $I$ is an ideal of $A$.

\medskip

\textbf{Definition.} An algebra $A$ is called \textit{Novikov algebra} if it satisfies the following identities:
\[(a,b,c)=(b,a,c), \eqno(1)\]
\[(ab)c=(ac)b. \eqno(2)\]

One can check that any Novikov algebra satisfies the following identities (see, for example, \cite{TUZ}):
\[(ad,b,c)=(a,bd,c)=(a,b,c)d.\eqno(3)\]

%\[(ad,b,c)=((ad)b)c-(ad)(bc)=((ab)c)d-(a(bc))d=(a,b,c)d.\]
%\[(a,bd,c)=(bd,a,c)=(b,a,c)d=(a,b,c)d.\]

\medskip

\textbf{Definition.} For any algebra some subspaces are defined. \textit{The nucleus}:
\[N(A)=\{n\in A\,|\, (n,a,b)=(a,n,b)=(a,b,n)=0 \,\,\forall a,b\in A\}, \eqno(4)\]
\textit{commutative center}:
\[K(A)=\{k\in A\,|\, [k,a]=0\,\forall a\in A\},\eqno(5)\]
and \textit{center}:
\[Z(A)=N(A)\cap K(A).\]
It is well known (\cite{Zhevl}) that nucleus and center are always subalgebras.\medskip

\textbf{Lemma 1.} \textit{Let $A$ be a Novikov algebra, $x,y,z\in A$ and $n\in N(A)$. Then 
\[n(x,y,z)=(nx,y,z)=(xn,y,z)=(x,yn,z)=\]
\[=(x,ny,z)=(x,y,nz)=(x,y,zn)=(x,y,z)n=0 \eqno(6)\]
}

\textbf{Proof.} Two equalities are true for any nonassociative algebra. First:
\[(x,y,zn)=(xy)(zn)-x(y(zn))\overset{(4)}{=}((xy)z)n-x((yz)n)\overset{(4)}{=}\]
\[\overset{(4)}{=}((xy)z)n-(x(yz))n=(x,y,z)n,\]
And second:
\[(nx,y,z)=((nx)y)z-(nx)(yz)\overset{(4)}{=}n((xy)z)-n(x(yz))=n(x,y,z).\eqno(7)\]
Two more are true for right-commutative algebras. First:
\[(x,y,nz)=(xy)(nz)-x(y(nz))\overset{(4)}{=}((xy)n)z-x((yn)z)\overset{(2)}{=}\]
\[\overset{(2)}{=}((xy)z)n-x((yz)n)
\overset{(4)}{=}((xy)z)n-(x(yz))n=(x,y,z)n,\]
And second:
\[(x,ny,z)=(x(ny))z-x((ny)z)\overset{(4)}{=}((xn)y)z-x(n(yz))\overset{(2)}{=}\]\[\overset{(2)}{=}((xy)z)n-x(n(yz))\overset{(4)}{=}((xy)z)n-(xn)(yz)\overset{(2)}{=}\]
\[\overset{(2)}{=}((xy)z)n-(x(yz))n=(x,y,z)n. \eqno(8)\]
We have two more identities for Novikov algebras from (3) (one of them is true even for right-commutative algebras):
\[(xn,y,z)\overset{(1)}{=}(x,yn,z)=(x,y,z)n,\]
But then for Novikov algebras we have
\[n(x,y,z)\overset{(7)}{=}(nx,y,z)\overset{(1)}{=}(y,nx,z)\overset{(8)}{=}(y,x,z)n\overset{(1)}{=}(x,y,z)n.\]
Also by (3)
\[(nx,y,z)\overset{(3)}{=}(n,y,z)x\overset{(4)}{=}0.\]
\qed
\medskip

\textbf{Definition.} If $A$ is an algebra then its \textit{associator ideal} $D(A)$ is an ideal, generated by all associators $(a,b,c)$, $a,b,c\in A$.  In \cite{TUZ} it has been shown that if $A$ is a Novikov algebra then $D(A)$ is spanned by the associators. Lemma 1 implies that $D(A)N(A)=N(A)D(A)=0$. 

\medskip

\textbf{Definition.} Let $A$ be an algebra and $M\subset A$. We will use the following notations for \textit{left annihilator} and \textit{right annihilator} of $M$ in $A$:
\[Ann_l(M)=\{a\in A\,|\, aM=0\}, \eqno(9)\]
\[Ann_r(M)=\{a\in A\,|\, Ma=0\}. \eqno(10)\]
\smallskip
\vspace{-10pt}

\textbf{Lemma 2.} \textit{Let $I$ be a left ideal in Novikov algebra $A$. Then $Ann_l(I)$ is an ideal of $A$. If $I$ is an ideal of $A$ then $Ann_r(I)$ is a left ideal of $A$.}

\textbf{Proof.} Let $x\in Ann_l(I)$, $a\in A$, $i\in I$. Then we have 
\[(xa)i\overset{(2)}{=}(xi)a\overset{(9)}{=}0,\]
\[(ax)i=(a,x,i)+a(xi)\overset{(1),(9)}{=}(x,a,i)=(xa)i-x(ai)\overset{(9)}{=}0.\]
Let $I$ be an ideal of $A$. If $y\in Ann_r(I)$ then 
\[i(ay)=-(i,a,y)+(ia)y\overset{(1),(10)}{=}-(a,i,y)=-(ai)y+a(iy)\overset{(10)}{=}0.\]
\qed

\medskip

\textbf{Lemma 3.} \textit{If $A$ is a Novikov algebra then $N(A)$ and $Z(A)$ are ideals of $A$. Also $K(A)=Z(A)$.}

\textbf{Proof.} Let $n\in N(A)$, $a\in A$. Then for any $b,c\in A$ we have by (6)
\[(na,b,c)=(a,nb,c)=(a,b,nc)=\dots =(n,b,c)a=0\]
So, $na,an\in N(A)$ and $N(A)$ is an ideal of $A$. Let $z\in Z(A)$. Then $za=az\in N(A)$ and
\[(za)x\overset{(4)}{=}(zx)a\overset{(5)}{=}(xz)a\overset{(4)}{=}x(za).\]
So, $za\in Z(A)$ and $Z(A)$ is an ideal of $A$.

Let $k\in K(A)$. Then
\[(x,k,y)\overset{(1)}{=}(k,x,y)\overset{(5)}{=}(xk)y-k(xy)\overset{(1)}{=}(xy)k-k(xy)\overset{(5)}{=}0,\]
\[(x,y,k)=(xy)k-x(yk)\overset{(1),(5)}{=}(xk)y-x(ky)=(x,k,y)=0.\]
So, $k\in N(A)\cap K(A)=Z(A)$. \qed

\medskip

\textbf{Definition.} An algebra $A$ is called \textit{prime} if for any ideals $I,J\unlhd A$ we have 
\[IJ=0 \to I=0 \text{ or } J=0.\]
Let us notice that product of ideals in Novikov algebra is an ideal.

\medskip

\textbf{Theorem 1.} \textit{If $A$ is a nonassociative prime Novikov algebra then \linebreak $Z(A)=N(A)=0$.}

\textbf{Proof.} Indeed, let $(a,b,c)\neq 0$. Then $D(A)\neq 0$ and $D(A)N(A)=0$. So, $N(A)=0$ by Lemma 3. \qed

\medskip

%\textbf{Следствие 1.} \textit{Пусть $A$ --- полупервичная неассоциативная алгебра Новикова. Тогда \linebreak $Z(A)=N(A)$.}

%\textbf{Доказательство.} Пусть $n\in N(A)$ и $A$ является подпрямым произведением первичных алгебр $A_i$. Тогда $(n,a,b)=0$ для любых $a,b\in A$ и $(\overline{n},\overline{a},\overline{b})=\overline{0}$, где $\overline{x}$ --- гомоморфный образ элемента $x$ в первичной алгебре $A_i$. Таким образом, $\overline{n}\in N(A_i)$. Если $A_i$ неассоциативна, то по теореме 1 $\overline{n}=\overline{0}$. Допустим, $A_i$ ассоциативна, тогда она является областью целостности. Но тогда $\overline{[n,a]}=[\overline{n},\overline{a}]=\overline{0}$ для любого $a\in A$ в любом гомоморфном образе $A_i$. Но тогда $[n,a]=0$ для любого $a\in A$, так что $n\in Z(A)$. Следствие доказано.

%\medskip

\section{Ideals of semiprime Novikov algebras.}

\textbf{Definition.} An algebra $A$ is called \textit{semiprime} if for any ideal $I$ we have
\[I^2=0 \to I=0.\]
An ideal $I$ is called \textit{trivial} if $I^2=0$. So, an algebra is semiprime iff it does not contain trivial nonzero ideals.

\medskip

\textbf{Lemma 4.} \textit{Let $A$ be a semiprime Novikov algebra, $I\unlhd A$, $V\unlhd I$ and $V^2=0$. Then $AV+V$ and $VA+V$ are trivial ideals of an algebra $I$.}

\textbf{Proof.} We have 
\[AI+IA\subset I, \qquad VI+IV\subset I \eqno(11)\]
and
\[V^2=0, \qquad (VI)V=(V,V,I)=(V,I,V)=(I,V,V)=0. \eqno(12)\]
Then we have
\[(AV)I\overset{(1)}{\subset} (AI )V\overset{(11)}{\subset} I V\overset{(11)}{\subset} V\]
\[I (AV)\subset (I ,A,V)+(I A)V\overset{(1),(11)}\subset (A,I ,V)+V\overset{(11)}{\subset} AV+V\]
For $VA$ we have
\[I (VA)\subset (I ,V,A)+(I V)A\overset{(1)}{\subset} (V,I,A)+(IV)A\overset{(11)}{\subset} \]
\[\overset{(11)}{\subset} (VI )A+V(I A)+VA\overset{(11)}{\subset} VA+V\]
\[(VA)I \overset{(1)}{\subset} (VI )A\overset{(11)}{\subset} VA.\]
So, $AV+V$ and $VA+V$ are ideals of $I$. 

\[(A,V,V)I\overset{(3)}{\subset}(AI,V,V)\overset{(11)}{\subset} (IV)V+IV^2\overset{(12)}{=}0,\]
i.e. $(a,u,v)\in Ann_l(I)$ for any $a\in A$, $u,v\in V$.
But $(a,u,v)\overset{(1)}{=}(u,a,v)\in V\subset I$, so
\[(a,u,v)\in I\cap Ann_l(I).\]
But $(I\cap Ann_l(I))^2=0$, so by Lemma~2 $I\cap Ann_l(I)=0$, hence $(a,u,v)=0$. So, we have
\[0=(a,u,v)=(au)v-a(uv)\overset{(12)}{=}(au)v,\]
i.e.
\[(AV)V=0\]
For $VA$ we have 
\[(VA)V\overset{(11)}{\subset} V^2A\overset{(12)}{=}0.\]
Then
\[(V,A,V)I\overset{(3)}{\subset}(V,AI,V)\overset{(11)}{\subset} (V,I,V)\overset{(11)}{\subset} V^2\overset{(12)}{=}0\]
$(V,A,V)\subset V\subset I$, so similarly we have $(u,a,v)=0$ for any $a\in A$, $u,v\in V$.

But 
\[0=(u,a,v)=(ua)v-u(av)\overset{(12)}{=}-u(av),\]
so 
\[V(AV)=0.\eqno(13)\]
\medskip

Then we have
\[(V(VA))I\subset (V,VA,I)+V((VA)I)\overset{(3)}{\subset}\] \[\overset{(3)}{\subset}(V,V,I)A+V(V,A,I)+V(V(AI)) \overset{(12)}{=}V(V,A,I) \overset{(1)}{\subset} V(A,V,I) =\]
\[=V((AV)I)+V(A(VI))\overset{(2),(11)}{\subset} V((AI)V)+V(AV)\overset{(13)}{=}0\]
Then 
\[V(VA)\subset Ann_l(I)=0 \eqno(14)\] 
So,
\[(VA)(VA)\overset{(2)}{\subset} (V(VA))A\overset{(14)}{=}0,\]
Now, we need to check that $(AV)(AV)=0$. We have
\[(A,V,AV)I\overset{(3)}{\subset}(AI,V,AV)\subset ((AI)V)(AV)+(AI)(V(AV))\overset{(13)}{=}0,\]
so
\[(A,V,AV)=0.\]
Then $(a,v,bu)=0$ for any $a,b\in A$, $v,u\in V$. Then we have for any $i\in I$
\[(a,v,bu)i\overset{(3)}{=}(ai,v,bu)=((ai)v)(bu)-(ai)(v(bu))\overset{(13)}{=}0,\]
so
\[(A,V,AV)I=0,\]
i.e. $(A,V,AV)\subset Ann_l(I)=0$ and
\[0=(a,u,bv)=(au)(bv)-a(u(bv))=\overset{(13)}{=}(au)(bv)\]
and $(AV)(AV)=0$. 
\qed

\bigskip

\textbf{Theorem 2.} \textit{Let $A$ be a semiprime Novikov algebra. Then its ideals are semiprime Novikov algebras.}

\textbf{Proof.} Let $I\unlhd A$. If $U\unlhd I$ is a trivial ideal of $I$ then by Zorn's lemma there exists maximal trivial ideal $V$ of an algebra $I$. But $VA+V$ is trivial too, hence $VA\subset V$. An ideal $AV+V$ is trivial too hence $AV\subset V$. It means that $V$ is a trivial ideal of $A$ hence $V=0$. So, an algebra $I$ is semiprime. \qed

\medskip

This theorem allow us to construct the first radical in the variety of Novikov algebras. 

\medskip

\textbf{Corollary 1.} \textit{Baer radical exists in the variety of Novikov algebras.}

\textbf{Proof.} Result is implied from Theorem 2 and Proposition 8.5. of \cite{Zhevl}. \qed

\bigskip

\section{Ideals of prime Novikov algebras.}

The following lemma have the same structure as Lemma 8.3. in \cite{Zhevl}, but the proof depends on identities of Novikov algebras.\smallskip

\textbf{Lemma 5.} \textit{Let $A$ be a Novikov algebra, $I$ be an ideal of $A$ and $M$ be an ideal of $I$. For any $a\in A$ we have:\\
a) $Ma+M$ and $aM+M$ are ideals of $I$;\\
b) $(MA)I^2\subset M$;\\
c) $(Ma)^2(Ma)^2\subset M$;\\
d) $(AM)I^2\subset M$;\\
e) $(aM)^2(aM)^2\subset M$.
}

\textbf{Proof.} We will use $MI+IM\subset M$ and $IA+AI\subset I$ without any clarifications. 

a) 
\[(Ma)I\overset{(2)}{\subset}(MI)a\subset Ma,\]
\[I(Ma)\subset (I,M,a)+(IM)a\overset{(1)}{\subset} (M,I,a)+Ma\subset Ma+M\]
So, $Ma+M$ is an ideal of $I$. 
\medskip

\[(aM)I\overset{(2)}{\subset}(aI)M\subset M\]
\[I(aM)\subset (I,a,M)+(Ia)M\overset{(1)}{\subset} (a,I,M)+M\subset aM+M.\]
So, $aM+M$ is an ideal of $I$.

\medskip

%г) \[(Ma)^2I\subset ((Ma)(Ma))I=((Ma)I)(Ma)=((MI)a)(Ma)\subset (Ma)^2\]

b) \[(MA)I^2\subset (MA,I,I)+((MA)I)I\overset{(3),(2)}{\subset} (M,IA,I)+((MI)A)I\subset\]
\[\subset (M,I,I)+(MI,A,I)+(MI)(AI)\overset{(3)}{\subset} M+(M,AI,I)\subset M \]

c) \[(Ma)^2(Ma)^2\overset{a)}{\subset} (M+Ma)(Ma)^2\subset MI^2+(Ma)I^2\overset{b)}{\subset} M.\]

d) \[(AM)I^2\overset{(2)}{\subset} (AI^2)M\subset M.\]

e) \[(aM)^2(aM)^2\overset{d)}{\subset} (M+aM)(aM)^2\subset MI^2+(aM)I^2\subset M+(aM)I^2\overset{d)}{\subset} M\]
\qed

\medskip

\textbf{Lemma 6.} \textit{Let $A$ be a Novikov algebra, $I$ be an ideal of $A$ and $M$ be an ideal of $I$. If an algebra $I/M$ is semiprime then $M$ is an ideal of $A$.}

\textbf{Proof.} Suppose that $M$ is not an ideal of $A$. Then $Ma$ or $aM$ does not contain in $M$ for some $a\in A$. Suppose $Ma$ does not contain in $M$. Then Lemma 5 implies that $M_1=M+Ma$ is an ideal in $I$. Denote its image in $I/M$ by $\overline{M_1}$. We have $((Ma)^2)^2\subset M$, so $(\overline{M_1}^2)^2=\overline{0}$. It means that $I/M$ contains trivial ideals. We have a contradiction, so $M$ is a right ideal in $A$. Similarly, $M$ is a left ideal in $A$. \qed

\medskip

\textbf{Lemma 7.} \textit{Let $I$ be a left ideal of a Novikov algebra $A$. Then a set $M=\{x\in I|xA\subset I\}$ is an ideal of an algebra $A$.}

\textbf{Proof.} Indeed, let $x\in M$, $a,b\in A$. Then
\[(xa)b=(x,a,b)+x(ab)\overset{(1)}{=}(a,x,b)+x(ab)=(ax)b-a(xb)+x(ab)\overset{(2)}{=}\]
\[\overset{(2)}{=}(ab)x-a(xb)+x(ab)\in I. \]
And obviously
\[(ax)b\overset{(2)}{=}(ab)x\in I.\]
So, $M$ is an ideal of $A$. \qed

\bigskip

\textbf{Theorem 3.} \textit{An ideal of a prime Novikov algebra is a prime Novikov algebra.}

\textbf{Proof.} Let $A$ be a prime Novikov algebra and let $I$ be an ideal of $A$. Theorem 2 implies that $I$ is semiprime. Suppose that $I$ contains two nonzero ideals $U,V\unlhd I$ with property $UV=0$. Then $(U\cap V)^2=0$ and $U\cap V=0$ by semiprimeness of $I$. It is enough to look a case $U=Ann_{l, I}(V)$. 

Let us denote $I/U=\overline{I}$. 

Let $\overline{T}$ is a trivial ideal of $\overline{I}$, $\overline{T}^2=\overline{0}$. If $T$ is a preimage of $\overline{T}$ then $T^2\subset U$. We have

\[(T\cap V)^2\subset T^2\cap V\subset U\cap V = 0.\]
So, $(T\cap V)^2=0$ and $T\cap V=0$. Then we have

\[VT\subset V\cap T = 0, \qquad TV\subset V\cap T=0 \to VT=TV=0.\]

We have $T\subset Ann_l(V)=U$ hence $\overline{T}=0$. So, the algebra $\overline{I}$ is semiprime. Lemma 6 implies that $U$ is an ideal of $A$ with nonzero right annihilator $Ann_r(U)$.  
\medskip

Let $M=\{x\in Ann_r(U)| xA\subset Ann_r(U)\}$ and $v\in V$. Then for any $u\in U$, $a\in A$ we have
\[u(va)=(u,v,a)-(uv)a\overset{(1)}{=}(v,u,a)=(vu)a-v(ua)=0.\]
So, $V\subset M$ and $M\neq 0$. Lemma 7 implies that $M$ is an ideal of $A$, so $UM$ is a nonzero product of ideals. It is a contradiction. \qed

\bigskip

\textbf{Theorem 4.} \textit{Let $A$ be a Novikov algebra and $I$ be a minimal ideal of $A$. Then either $I^2=0$, or $I$ is a simple algebra.}

\textbf{Proof.} We have $I^2\subset I$ and $I^2$ is an ideal of $I$. So either $I^2=0$, or $I^2=I$. We need to crack a case of $I^2=I$.

Let $B$ be an ideal of an algebra $I$ and $B\neq I$. Denote $C_1=(IB)I+I(IB)$. We want to prove that $C_1$ is an ideal of $A$. 
\[((IB)I)A\overset{(2)}{\subset} ((IA)B)I\subset (IB)I\subset C_1,\]
\[(I(IB))A\overset{(2)}{\subset} (IA)(IB)\subset I(IB)\subset C_1.\]
So, $C_1$ is a right ideal of $A$. 
\[A(I(IB))\subset (A,I,IB)+(AI)(IB)\overset{(1)}{\subset} (I,A,IB)+I(IB)=\] 
\[=(I^2,A,IB)+I(IB)\overset{(3)}{\subset} (I,AI,IB)+I(IB)\subset (I,I,B)+I(IB)\subset \]
\[\subset I^2B+I(IB)\overset{(2)}{\subset} (IB)I+I(IB)=C_1,\]
\[A((IB)I)\subset (A,IB,I)+(A(IB))I\overset{(1)}{\subset}\] \[\overset{(1)}{\subset}(IB,A,I)+(A,I,B)I+((AI)B)I\subset\] \[\subset ((IB)A)I+(IB)(AI)+(A,I,B)I+(IB)I\overset{(2)}{\subset}\] \[\overset{(2)}{\subset} ((IA)B)I+(IB)I+(IB)I+(A(IB))I\overset{(2)}{\subset} (IB)I+(AI)(IB)\subset C_1.\]
So, $C_1$ is an ideal of $A$ and $C_1\subset B$. So, $C_1\subsetneq I$ and $C_1=0$ by minimality of $I$. Hence $(IB)I=0$. 

Suppose that $I\cap Ann_l(I)\neq 0$. But then $I\subset Ann_l(I)$ and $I^2=0$. A contradiction, so $I\cap Ann_l(I)=0$. 

But $IB\subset I\cap Ann_l(I)=0$ hence $IB=0$.

Denote $C_2=I(BI)$. We want to prove that $C_2$ is an ideal of $A$. 
\[(I(BI))A\overset{(2)}{\subset} (IA)(BI)\subset I(BI)= C_2,\]
\[A(I(BI))\subset (A,I,BI)+(AI)(BI)\overset{(1)}{\subset} (I,A,BI)+I(BI)=\]
\[=(I^2,A,BI)+I(BI)\overset{(3)}{\subset} (I,AI,BI)+I(BI)\subset \]
\[\subset I(BI)+I(I(BI))\subset 
I(BI)+I(I,B,I)\overset{(1)}{\subset}\]
\[\overset{(1)}{\subset} I(BI)+I(B,I,I)\subset I(BI)+I((BI)I)+I(BI^2)\subset\]
\[\subset I(BI)+(I,BI,I)+(I(BI))I\overset{(3)}{\subset} I(BI)+(I^2,B,I)\subset I(BI)=C_2.\]
So, $C_2$ is an ideal of $A$ and $C_2\subset B$. So, $C_2\subsetneq I$ and $C_2=I(BI)=0$ by minimality of $I$. Hence $(B,I,I)\overset{(1)}{\subset}(I,B,I)\subset(IB)I+I(BI)=0$.

Denote $C_3=(BI)I$. We want to prove that $C_3$ is an ideal of $A$. 
\[((BI)I)A\overset{(2)}{\subset}  (BI,A,I)+(BI)(AI)\overset{(3)}{\subset} (B,I,I)+(BI)I=C_3,\]
\[A((BI)I)\subset (A,BI,I)+(A(BI))I\overset{(3),(2)}{\subset}\] \[\overset{(3),(2)}{\subset}(I,B,I)+(AI)(BI)\subset I(BI)=0.\]
So, $C_2$ is an ideal of $A$ and $C_2\subset B$. So, $C_2\subsetneq I$ and $C_2=(BI)I=0$ by minimality of $I$. But then $BI\subset Ann_l(I)\cap I = 0$ hence $BI=0$. But then $B\subset Ann_l(I)\cap I=0$ hence $B=0$. \qed

\bigskip

Let us notice that commutative Novikov algebra is associative.

\medskip

\textbf{Lemma 8.} \textit{Let $A$ be a Novikov algebra with a nonzero commutative ideal $H$ and $H^2=H$. Then $H\subset N(A)$. }

\textbf{Proof.} $H$ is an associative and commutative algebra. Let $h,h_i\in H$, $a,b\in A$. Then $(h,a,b)=(a,h,b)$ is a linear combination of elements like
$(h_1h_2,a,b)=(h_1,ah_2,b)=(h_1,h_3,b).$
Then we have
\[(h_1,h_3,b)=(h_3h_1)b-h_1(h_3b)=(h_3b)h_1-h_1(h_3b)=0.\]
Last equation is true because $h_3b\in H$. For the third argument we have
\[(a,b,h_1h_2)=(ab)(h_1h_2)-a(b(h_1h_2))=\]
\[=-(ab,h_1,h_2)+((ab)h_1)h_2+a(b,h_1,h_2)-a((bh_1)h_2)=\]
\[=((ah_1)b)h_2-a((bh_1)h_2)=(ah_1)(bh_2)-a((bh_2)h_1)=\]
\[=(a(bh_2))h_1-a((bh_2)h_1)=0.\]
So, $H\subset N(A)$. \qed

\bigskip

\textbf{Corollary 2.} \textit{If $A$ is a prime nonassociative Novikov algebra then it does not contain nonzero commutative ideals with property $H=H^2$.}

\textbf{Proof.}  Theorem 1 implies $N(A)=0$. Then Lemma~8 implies $H=0$. \qed

\bigskip

\textbf{Corollary 3.} \textit{Prime Novikov algebra with a minimal commutative ideal is associative.}

\textbf{Proof.} If $A$ contains minimal commutative ideal $H$ then $H^2=H$ by semiprimeness. Corollary 2 implies the result. \qed

\end{document}